\def \version {17 October 2015}

\documentclass[12pt]{article}
\usepackage{amssymb}

\usepackage{tikz,ifthen}

\newcommand{\qed}{\hfill ~$\square$}

\def \nsp {\vspace{-.5ex}}
\def \nspp {\vspace{-2ex}}

\newtheorem{Theorem}{Theorem}
\newtheorem{lem}[Theorem]{Lemma}
\newtheorem{defi}[Theorem]{Definition}
\newtheorem{crl}[Theorem]{Corollary}
\newtheorem{prp}[Theorem]{Proposition}
\newtheorem{prm}[Theorem]{Problem}
\newtheorem{rmk}[Theorem]{Remark}
\newtheorem{xmp}[Theorem]{Example}
\newtheorem{conjecture}[Theorem]{Conjecture}

\def \taut {\tau_\vartriangle}
\def \nut {\nu_\vartriangle}

\def \bp {\begin{prp} \ }
\def \ep {\end{prp}}

\def \bpm {\begin{prm} \ }
\def \epm {\end{prm}}

\def \bc {\begin{crl} \ }
\def \ec {\end{crl}}

\def \thm {\begin{Theorem} \ }
\def \ethm {\end{Theorem}}

\def \bl {\begin{lem} \ }
\def \el {\end{lem}}

\def \bd {\begin{defi} \ \rm }
\def \ed {\end{defi}}

\def \brm {\begin{rmk} \ }
\def \erm {\end{rmk}}

\def \bxm {\begin{xmp} \ \rm }
\def \exm {\end{xmp}}

\def \nmr {\begin{enumerate}}
\def \enmr {\end{enumerate}}

\def \tmz {\begin{itemize}}
\def \etmz {\end{itemize}}

\def \nin {\noindent}
\def \bsk {\bigskip}
\def \msk {\medskip}

\def \pf {\nin{\bf Proof.} \ }
\def \prf {\nin{\it Proof.} \ }

\def \NP {{\sf NP}}

\def \taut {\tau_\vartriangle}
\def \nut {\nu_\vartriangle}

\def\cF{{\cal F}}

\def\cT{{\cal T}}

\newtheorem{claim}{Claim}

\def \bcl {\begin{claim} \ }
\def \ecl {\end{claim}}

\newtheorem{con}{Condition}

\def \bcon {\begin{con} \ \rm }
\def \econ {\end{con}}

\def\T{\mathcal{T}}
\def \dia {\hfill $\Diamond$}

\author{$^1$Aparna Lakshmanan S, \, $^{2,3}$\,Csilla Bujt\'{a}s, \,$^{2,3}\,$Zsolt Tuza, \,\\
\\
\normalsize $^1$\,Department of \nsp Mathematics\\
\normalsize St.~Xavier's College for \nsp Women, Aluva, India\\
\small {\tt E-mail: \nspp aparnaren@gmail.com}\\
\\
\normalsize $^2$\,Department of Computer Science and Systems
  \nsp Technology\\
\normalsize University of \nsp Pannonia, Veszpr\'em, Hungary \\
\small {\tt E-mail: \nspp bujtas@dcs.uni-pannon.hu}\\
\\
\normalsize $^3$\,Alfr\'ed R\'enyi \nsp Institute of Mathematics\\
\normalsize Hungarian Academy of \nsp Sciences, Budapest, Hungary \\
\small {\tt E-mail: tuza@dcs.uni-pannon.hu}\\
\\}

\date{~ \vspace{-6ex} ~}

\title{%{~ \vspace{-2ex} ~} \\
 Induced cycles in triangle graphs%~\thanks{
 }

\date{\small Latest update \vspace{-2ex} on \version }

\begin{document}

\maketitle

\begin{abstract}
The triangle graph of a graph $G$, denoted by $\T(G)$, is the graph
whose vertices represent the triangles
 ($K_3$ subgraphs) of $G$, and two vertices of
$\mathcal{T}(G)$ are adjacent if and only if the corresponding
triangles share an edge. In this paper, we characterize graphs
 whose triangle graph is a cycle and then extend the result to
obtain a  characterization of $C_n$-free triangle graphs. As a
consequence, we give a forbidden subgraph characterization of
graphs $G$ for which $\T(G)$ is a tree, a chordal graph, or a
perfect graph. For the class of graphs whose triangle graph is
perfect, we verify a conjecture of the third author concerning
packing and covering of triangles.

\bsk

\nin {\bf Keywords :} Triangle graph, $F$-free graph, perfect
graph, Tuza's Conjecture

\bsk

\nin AMS Subject Classification:
  {%\bf
   05C76,
%% (Graph operations (line graphs, products, etc.)),
  05C75
%% (Structural characterization of types of graphs),
 (primary);
%%  05C62 Graph representations (geometric and intersection representations, etc.)
  05C70
%% (Factorization, matching, partitioning, covering and packing)
 (secondary).
  }

\end{abstract}

\section{Introduction}

In a simple undirected graph $G$, a {\it triangle\/} is a complete
subgraph on three vertices. The {\it triangle graph\/} of $G$,
denoted by $\T(G)$, is the graph whose vertices represent the
triangles of $G$, and two vertices of $\mathcal{T}(G)$ are adjacent
if and only if the corresponding triangles of $G$ share an edge.
This notion was introduced independently several times under
different names and in different contexts  \cite{Pul, T-Oberw, Ega,
Bal}.
One fundamental motivation is its obvious relation to the important
 class of line graphs.

\bsk

In a more general setting, for a $k \ge 1$, the $k$-line graph
$L_k(G)$ of $G$ is a graph which has vertices corresponding to the
$K_k$ subgraphs of $G$, and two vertices are adjacent in $L_k(G)$ if
the represented $K_k$ subgraphs of $G$ have $k-1$ vertices in
common. Hence,  2-line graph means line graph in the usual sense,
whilst 3-line graph is just the triangle graph, which is our current
subject.

\bsk

Beineke's classic result \cite{Bei} gave a characterization of $2$-line
 graphs in terms of nine forbidden subgraphs. This implies that $2$-line graphs
  can be recognized in polynomial time.
  In contrast to this, as  proved very recently in \cite{ABT2},  the recognition problem of triangle
 graphs (and also, that of $k$-line graphs for each  $k \ge 3$)
 is \NP-complete. In the same paper \cite{ABT2}, a necessary and
 sufficient condition is given for nontrivial connected graphs
  $G$ and $H$ to ensure that their
 Cartesian product $G\Box H$ is a triangle graph.

\bsk

 Further related results have
% recently
  been obtained
 by Laskar, Mulder and Novick  \cite{Laskar1}. They prove that for an `edge-triangular'
  and `path-neighborhood' graph $G$
 (that is when the open neighborhood of $v$ induces a non-trivial path for  each
 vertex $v \in V(G)$), the triangle graph $\cT(G)$ is a tree if and only if $G$ is
 maximal outerplanar. Also, they raise the characterization problem
 of a path-neighborhood graph $G$ for which $\cT(G)$ is a cycle (\cite[Problem 3]{Laskar1}).
 As an immediate consequence of our Theorem~\ref{cycle}, we will answer this
 question; moreover we will give a forbidden subgraph
 characterization of graphs whose triangle graph is a tree.

 Triangle graphs were studied from several further aspects;
  see e.g.\ \cite{Bag, Bal,  Ega, Le, Eri, Pris, Ram1, Ram2, Ram3}.\\

\subsection{Standard definitions}

 Given a graph $F$, a graph $G$ is called {\it $F$-free\/} if no {\it
induced\/}  subgraph of $G$ is isomorphic to $F$. When $\cF$ is a
set of graphs, $G$ is {\it $\cF$-free\/} if it is $F$-free for all
$F\in \cF$.
 On the other hand, when we say that a graph $F$ is a
\emph{forbidden subgraph}
 for a class $\cal G$ of graphs, it means that no $G\in {\cal G}$ may
 contain \emph{any} subgraph isomorphic to $F$. \\

As usual, the complement of a graph $G$ is denoted by $\overline{G}$.
 The $n$th power    of a graph $G$ is the
graph $G^n$    whose vertex set is $V(G^n)=V(G)$ and two vertices
are adjacent in $G^n$ if and only if their distance is at most $n$
in $G$. Moreover, given two graphs $G_1=(V_1,E_1)$ and
$G_2=(V_2,E_2)$,
 we use the notation $G_1 \vee G_2$ for
 the \emph{join} of $G_1$ and $G_2$, that is a graph with one copy of
$G_1$ and $G_2$ each, being vertex-disjoint,
 and all the vertices of $G_1$ are made
 adjacent with all the vertices of $G_2$.
In particular, the {\it
$n$-wheel\/} $W_n$ ($n\geq 3$) is a graph
 $K_1 \vee C_n$ (where, as usual, $K_n$ and $C_n$ denote the
 $n$-vertex complete graph and the $n$-cycle, respectively).
An {\it odd wheel\/} is a graph $W_n$ where $n\geq3$ is odd;
 and an {\it odd hole\/} in a graph is an
{\it induced\/} $n$-cycle of odd length $n\ge 5$, whereas an {\it
odd anti-hole\/} is the complement of an odd hole.\\

While an {\it acyclic\/} graph does not contain any cycles,
 a {\it chordal graph\/} is a graph which
 does not contain \emph{induced} $n$-cycles for $n\geq4$.
The {\it chromatic number\/} $\chi(G)$ of a
graph $G$ is the minimum number of colors required to color the
vertices of $G$ in such a way that no two adjacent vertices
receive the same color. A set of vertices is {\it
independent\/} if all pairs of its vertices are non-adjacent. The
{\it independence number\/} $\alpha(G)$ of $G$ is the maximum
cardinality of an independent vertex set in $G$. A {\it clique\/}
is a complete subgraph maximal under inclusion (i.e., in our
 terminology different cliques in the same graph may have
 different size). The {\it clique
number\/} $\omega(G)$ is the maximum number of vertices of a
clique in $G$. The {\it clique covering number\/} $\theta(G)$ is
the minimum cardinality of a set of cliques that covers all
vertices of $G$. A graph $G$ is {\it perfect\/} if $\chi(G') =
\omega(G')$ for every induced subgraph $G'$ of $G$.

As usual, the {\it open neighborhood} $N(v)$ of $v$ is the set of
neighbors of $v$, whilst its {\it closed neighborhood} is $N[v]=N(v)
\cup \{v\}$. In a less usual way, we also refer to the subgraphs
induced by them as $N(v)$ and $N[v]$, respectively.

 Throughout this paper, the notation $K_n-G$ will refer
to the graph obtained from the complete graph $K_n$ by deleting the
\emph{edge set} of a subgraph isomorphic to $G$. In this way, for instance,
$K_4-K_3$ means the claw $K_{1,3}$.

\subsection{New definitions and terminology}   \label{defterm}

In this paper, we use the following special terminology for some
types of graphs.

\tmz
  \item The \emph{elementary types} are:
  \tmz
      \item[$(a)$] the wheel $W_4$,
      \item[$(b)$] the square $C_n^2$ of a cycle of length $n\geq 7$.
  \etmz

  \item  The {\it supplementary types} are the following ones. (For illustration, see
  Fig.~\ref{fig1}.)
\etmz \tmz
  \item[]$(A)$ \ $S_A=(V_A,E_A)$, where
  $V_A=\{v_i,u_i\mid 1\leq i\leq 4\}$    and
  $$E_A= \{v_iv_{i+1}\mid 1\leq i\leq  4\}\cup
  \{u_iv_{i-1},u_iv_{i }, u_iv_{i+1}\mid 1\leq i\leq  4\}$$
  \qquad \ (subscript addition taken modulo 4).

  \item[]$(B)$ \ $S_B=(V_B,E_B)$, where
  $V_B=\{v_i \mid 1\leq i\leq
  5\}\cup \{u_1,u_2,u_3\}$    and
  $$E_B= \{v_iv_{i+1}\mid 1\leq i\leq  5\}
   \cup \{v_3v_5,v_4v_1\} \cup
  \{u_iv_{i-1},u_iv_{i }, u_iv_{i+1}\mid 1\leq i\leq  3\}$$
  \qquad \ (subscript addition taken modulo 5).

  \item[]$(C)$ \ $S_C=(V_C,E_C)$, where
  $V_C=\{v_i \mid 1\leq i\leq
  6\}\cup \{u_1,u_2\}$    and
  $$E_C= \{v_iv_{i+1}\mid 1\leq i\leq  6\}
    \cup \{v_2v_4,v_3v_5,v_4v_6,v_5v_1\} \cup
  \{u_iv_{i-1},u_iv_{i }, u_iv_{i+1}\mid i=1,2\}$$
  \qquad \ (subscript addition taken modulo 6).

  \item[]$(D)$ \ $S_D=(V_D,E_D)$, where
  $V_D=\{v_i \mid 1\leq i\leq
  6\}\cup \{u_1,u_4\}$    and
  $$E_D= \{v_iv_{i+1}\mid 1\leq i\leq  6\}
   \cup \{v_1v_3,v_2v_4,v_4v_6, v_5v_1\} \cup
  \{u_iv_{i-1},u_iv_{i }, u_iv_{i+1}\mid i=1,4\}$$
  \qquad \ (subscript addition taken modulo 6).
\etmz

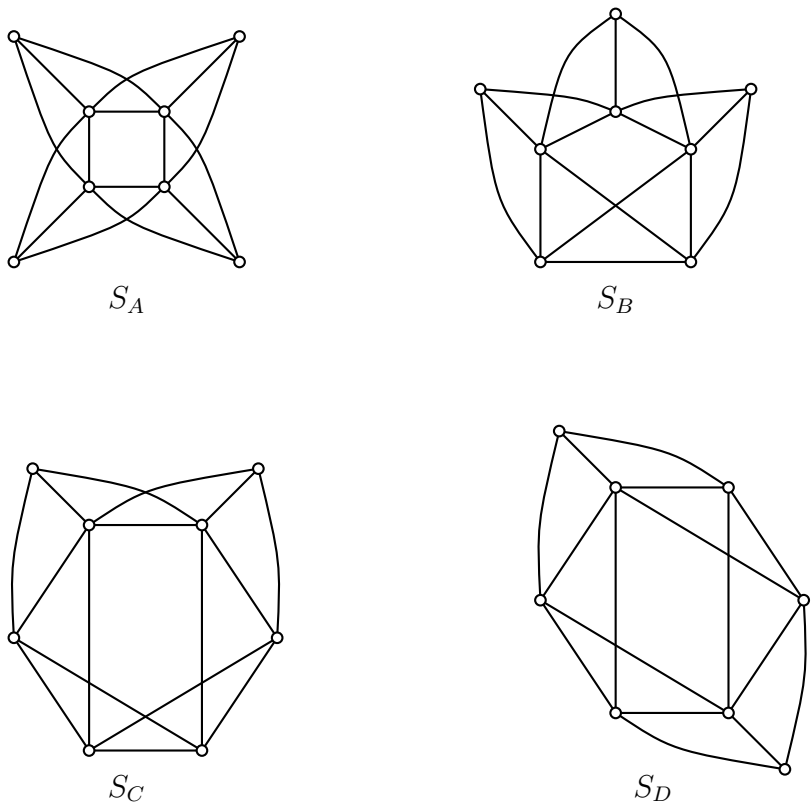
\begin{figure}[ht!] \label{fig1}
\begin{center}
\begin{tikzpicture}[scale=1,style=thick]
\def\vr{2pt} % \vr = vertex radius;  Set \vr = 2/scale for uniform sizing of vertices

%% vertices defined %%
 \path (0,10) coordinate (a0);
 \path (3,10) coordinate (a1);
  \path (1,11) coordinate (a2);
  \path (2,11) coordinate (a3);
  \path (1,12) coordinate (a4);
  \path (2,12) coordinate (a5);
  \path (0,13) coordinate (a6);
   \path (3,13) coordinate (a7);

   \path (7,10) coordinate (b5);
   \path (9,10) coordinate (b4);
   \path (7,11.5) coordinate (b1);
   \path (9,11.5) coordinate (b3);
   \path (8,12) coordinate (b2);
   \path (6.2,12.3) coordinate (b0);
   \path (9.8,12.3) coordinate (b7);
   \path (8,13.3) coordinate (b6);

   \path (0,5) coordinate (c0);
   \path (1,6.5) coordinate (c1);
   \path (2.5,6.5) coordinate (c2);
   \path (3.5,5) coordinate (c3);
   \path (2.5,3.5) coordinate (c4);
   \path (1,3.5) coordinate (c5);
   \path (0.25,7.25) coordinate (c6);
   \path (3.25,7.25) coordinate (c7);

   \path (7,5.5) coordinate (d0);
   \path (8,7) coordinate (d1);
   \path (9.5,7) coordinate (d2);
   \path (10.5,5.5) coordinate (d3);
   \path (9.5,4) coordinate (d4);
   \path (8,4) coordinate (d5);
   \path (7.25,7.75) coordinate (d6);
   \path (10.25,3.25) coordinate (d7);

%
%% edges %%
 \draw (a0) -- (a2) -- (a3) -- (a1);
 \draw (a4) -- (a2);
 \draw (a3) -- (a5) -- (a7);
 \draw (a5) -- (a4) -- (a6);
 \draw (b5) -- (b1) -- (b2) -- (b3) -- (b4) -- (b5);
 \draw (b0) -- (b1)-- (b4);
 \draw (b6) -- (b2);
 \draw (b7) -- (b3)-- (b5);

 \draw (c0) -- (c4)-- (c2) -- (c1) -- (c0) -- (c5) -- (c4) -- (c3) -- (c2) --(c7);
 \draw (c6) -- (c1)-- (c5) -- (c3);
 \draw (d0) -- (d4)-- (d2) -- (d1) -- (d0) -- (d5) -- (d4) -- (d3) -- (d2);
 \draw (d5) -- (d1)-- (d3);
 \draw (d6) -- (d1);
 \draw (d7) -- (d4);

 \draw (a0) .. controls (1.5,10.5)  .. (a3);
 \draw (a0) .. controls (0.5,11.5)  .. (a4);
 \draw (a1) .. controls (1.5,10.5)  .. (a2);
 \draw (a1) .. controls (2.5,11.5)  .. (a5);
 \draw (a7) .. controls (1.5,12.5)  .. (a4);
 \draw (a7) .. controls (2.5,11.5)  .. (a3);
 \draw (a6) .. controls (1.5,12.5)  .. (a5);
 \draw (a6) .. controls (0.5,11.5)  .. (a2);
 %%%%%%%%%%%%%%%%%%%%%%%%
 \draw (b0) .. controls (6.4,10.8)  .. (b5);
 \draw (b0) .. controls (7.5,12.2)  .. (b2);
 \draw (b7) .. controls (8.5,12.2)  .. (b2);
 \draw (b7) .. controls (9.6,10.8)  .. (b4);
 \draw (b6) .. controls (7.3,12.8)  .. (b1);
 \draw (b6) .. controls (8.7,12.8)  .. (b3);
 \draw (c6) .. controls (-0.05,6)  .. (c0);
 \draw (c6) .. controls (1.75,7)  .. (c2);
 \draw (c7) .. controls (3.55,6)  .. (c3);
 \draw (c7) .. controls (1.75,7)  .. (c1);
\draw (d6) .. controls (6.95,6.5)  .. (d0);
 \draw (d6) .. controls (8.75,7.5)  .. (d2);
 \draw (d7) .. controls (10.55,4.5)  .. (d3);
 \draw (d7) .. controls (8.75,3.5)  .. (d5);
%\draw (a0) .. controls (1,-1) and (4,-1) .. (a5);
%
%% vertices %%%
 \foreach \i in {0,1,...,7} {
     \draw (a\i)  [fill=white] circle (\vr);
  }
  \draw  (1.5,9.5) node  {$S_A$};
  \foreach \i in {0,1,...,7} {
     \draw (b\i)  [fill=white] circle (\vr);
  }
  \draw  (8,9.5) node  {$S_B$};
%% text %%%%%%%%%%%%%%%%%%%%%%%%%%%%%%%%%

\foreach \i in {0,1,...,7} {
     \draw (c\i)  [fill=white] circle (\vr);
     \draw (d\i)  [fill=white] circle (\vr);
  }
  \draw  (1.5,3) node  {$S_C$};
  \draw  (8.5,3) node  {$S_D$};

%% text %%%%%%%%%%%%%%%%%%%%%%%%%%%%%%%%%
%\draw [left] (a6) node {$u$}; \draw [right] (a16) node {$w$};
%% vertices %%%
%\foreach \i in {0,1,2,3} {
 %   \draw (a\i)  [fill=white] circle (\vr);
  %  \draw (c\i)  [fill=white] circle (\vr);
  %}
   % \draw (b0)  [fill=white] circle (\vr);
    %\draw (b1)  [fill=white] circle (\vr);
    %\draw (d1)  [fill=white] circle (\vr);
    %\draw (d2)  [fill=white] circle (\vr);
    %\draw (d3)  [fill=white] circle (\vr);
\end{tikzpicture}
\end{center}
\caption{The four graphs of supplementary type} \label{fig1}
\end{figure}

We also define two operations as follows.
 \tmz
  \item Suppose that $e=xy$ is an edge \underline{contained in exactly one
  triangle} $xyz$, whilst $xz$ and $zy$ belong to more than one triangle.
  An \emph{edge splitting} of $e$ means
  replacing $e$ with the 3-path $xwy$ (where $w$ is a new vertex)
  and inserting the further edge $wz$.
  \item Let $u$ and $v$ be two vertices \underline{at distance at least
  4  apart}. The \emph{vertex sticking} of $u$ and $v$ means
  removing $u$ and $v$ and introducing a new vertex $w$
  adjacent to the entire $N(u)\cup N(v)$.\footnote{`Vertex sticking'
  and its inverse operation were also introduced in \cite{Laskar1}.}
  \etmz

 The inverses of these operations can also be introduced in a natural way.
 \tmz
  \item Suppose that $xwz$ and $ywz$ are two triangles in the
   following position:
    $w$ has degree 3, $z$ is the unique common neighbor of $w$ and $x$
     and also of $w$ and $y$ (in particular, $x$ and $y$ are
     not adjacent), and $w$ and $z$ are the only common
     neighbors of $x$ and $y$.
  The \emph{inverse edge splitting} at $w$ means deleting $w$ and its
   three incident edges, and inserting the new edge $xy$.
  \item Let $w$ be a vertex whose neighborhood $N(w)$ is disconnected.
   The \emph{inverse vertex sticking} at $w$ means deleting $w$ and its
    incident edges, and inserting two new vertices $u$ and $v$ in
    such a way that $N(u)\cup N(v)=N(w)$, $N(u)\cap N(v)=\emptyset$,
    and each  component of $N(w)$ is either inside $N(u)$
    or inside $N(v)$.
 \etmz

\subsection{Our results}

In Section 2, we characterize the graphs whose triangle graph is a
cycle (Theorem \ref{cycle}) and then conclude a characterization of
 path-neighborhood graphs $G$ with $\cT(G) \cong C_n$ for some $n$.
 The latter one (Corollary \ref{png-char}) solves a problem
raised in \cite{Laskar1}.

 In Section 3, we prove a forbidden subgraph characterization of
 graphs $G$ with $C_n$-free triangle graphs for any specified $n \ge 3$
 (Theorem \ref{main}). Applying this result, we give necessary and
sufficient conditions for graphs $G$ whose $\T(G)$ is a tree, a
chordal graph, and a perfect graph, respectively.
 In a sense these results form a hierarchy since every
 tree is chordal, and
  every chordal graph is perfect \cite{Dir}.

  In Subsection~\ref{subsec:3-3} we consider the following old conjecture
   (usually referred to as ``Tuza's Conjecture'')
 of the third author  regarding packing and covering
 the triangles of a graph. It was formulated in 1981 \cite{Tuz81}.

%\nin {\bf Conjecture {\cite{Tuz81}:}}\ \ {\it
\begin{conjecture} \label{T-conj}
 If a graph\/ $G=(V,E)$ does not contain more than\/ $t$ mutually
  edge-disjoint triangles, then there exists\/ $E'\subseteq E$
  such that\/ $|E'|\le 2t$ and each triangle of\/ $G$ has at least
  one edge in\/ $E'$.
  \end{conjecture}

We prove that Conjecture~\ref{T-conj} holds for  graphs whose
triangle graph is perfect (Theorem \ref{perf-taunu}).

\section{Graphs whose triangle graph is a cycle}

In this section, we give a characterization of graphs whose triangle
graph is a cycle. We assume that every edge of $G=(V,E)$ is
contained in a triangle, and there are no isolated vertices.
  Before stating the theorem, let us prove that the required
 property is invariant under the two operations introduced in
 Section \ref{defterm}.

\bl \label{operations}
  For a graph\/ $G$, let\/ $G'$ be a graph obtained from\/ $G$ by
splitting an edge or sticking two vertices. Then\/ $\T(G')$ is a
cycle if and only if\/ $\T(G)$ is a cycle. \el

\pf Let first $G'$ be the graph obtained from $G$ by splitting an
edge $e=uv$ to the path $uwv$, where $e$ belongs to exactly one
triangle $uvx$. Let $t_1,t_2,\dots,t_n$ be the triangles in $G$.
Without loss of generality we may assume that $t_1 = uvx$, $t_2$
contains the edge $vx$ and $t_n$ contains the edge $ux$. Hence, by
the edge splitting, two neighboring triangles $uwx$ and  $wvx$
arose, where $uwx$ and $t_n$,  moreover $wvx$ and $t_2$  also have a
common edge.
 Then, clearly,
$t_1,t_2,\dots,t_n$ is an $n$-cycle in $\T(G)$ if and only if the
triangles  $uwx, wvx, t_2, \dots,t_n$ of $G'$ form an
$(n+1)$-cycle in $\T(G')$   in this cyclic order. Therefore,
$\mathcal{T}(G) \cong C_n$ if and only if $\mathcal{T}(G') \cong
C_{n+1}$.

Let now $G'$ be the graph obtained from $G$ by sticking two vertices
$u$ and $v$ to a new vertex $w$, where the distance between $u$ and
$v$ is at least four. If a triangle $t$ of $G$ contains neither $u$
nor $v$, then $t$ is a triangle in $G'$. If $t$ is a triangle which
contains $u$, say $t=ux_1x_2$, then $wx_1x_2$ is a triangle in $G'$.
The same holds for triangles containing $v$. Since $u$ and $v$ are
at distance at least 4 apart, no triangle of $G$ is damaged and no
new triangle can arise when $u$ and $v$ are stuck. Moreover, two
triangles share an edge in $G$ if and only if the corresponding
triangles have a common edge in $G'$. Therefore, $\T(G)\cong
\T(G')$. \qed

\bsk

 Observe that splitting an edge increases the number of triangles by exactly one,
 whilst sticking two vertices far enough does not change the number of
 triangles and the number of edges in a graph but increases edge density.
These observations have the following simple but important consequence.

\begin{crl}[Finite Reduction Lemma]
 For every fixed\/ $n$ there is a finite\/ $s_n$ such that,
starting from\/ \emph{any} graph whose triangle graph is a cycle,
 after\/ $s_n$ applications of edge splitting and vertex sticking in any feasible
 order, the length of the cycle\/ $\T(G)$ of the graph\/ $G$ obtained
 surely exceeds\/ $n$.
Equivalently, starting from any\/ $G$ whose triangle graph is a cycle,
 inverse edge splitting and inverse vertex sticking can be applied
 only finitely many times.
\end{crl}

 If an edge $e$ belongs to exactly one triangle,
  we call $e$ a \emph{private edge} (of that triangle);
 and if $e$ is contained in exactly two
 triangles, it is a \emph{doubly covered edge}.

\thm \label{cycle}   Let\/ $G$ be a graph  which contains no
isolated vertices and whose every edge is contained in at least one
triangle.
 Then,\/ $\mathcal{T}(G) \cong C_n$ for some\/
$n \geq 3$ if and only if
  \tmz \item[$(i)$] $G\cong K_5-K_3$, or
  \item[$(ii)$] $G$ is one of the elementary types or supplementary
  types, or
  \item[$(iii)$] $G$ can be obtained from one of
the   elementary types or supplementary types by a sequence of edge
splittings and vertex stickings. \etmz
 Moreover, graphs whose
triangle graph is an odd hole are characterized
%in the same way,
by\/ $(ii)$ and\/ $(iii)$, with properly chosen parity of the number of
edge splittings.
   \ethm

\pf Clearly, $\T( K_5-K_3) \cong C_3$, $\mathcal{T}(W_4)\cong C_4$,
and $\mathcal{T}(C_n^2)\cong C_n$ if $n\geq 7$, moreover the
triangle graphs of $S_A$, $S_B$,
 $S_C$ and $S_D$ are   isomorphic to $C_8$. By
 Lemma~\ref{operations}, also the triangle graphs of
  graphs satisfying $(iii)$ are cycles.

  \msk

  Now, assume that for a graph
  $G=(V,E)$,  fulfilling the conditions of the theorem,
its   triangle  graph $\T(G)$ is a cycle $C_n$.
   If an edge $e\in E$ is contained in more than two
triangles, then those triangles induce a complete subgraph
  of order at least 3 in the
triangle graph. Hence, if the latter is a cycle, then there cannot
be more than three triangles, thus $\mathcal{T}(G) = C_3$ and
$G\cong K_5-K_3$.

On the other hand, if a triangle has no  private edge    then the
degree of the corresponding vertex in $\mathcal{T}(G)$ will be at
least three, which contradicts the assumption that $\mathcal{T}(G)$
is a cycle.

\msk  From now on,   assume that $G\ncong K_5-K_3$ and consequently,
each triangle has precisely one private edge and exactly two  doubly
covered edges. Moreover, we will suppose that the inverse operations
of edge splitting and vertex sticking cannot be applied to $G$.
  (Due to the Finite Reduction Lemma above, we may assume this without
 loss of generality.) It will
suffice to prove that such non-reducible graphs necessarily belong
to an elementary or a supplementary type.

\bsk

 \nin
{\it Claim 1. } For every vertex\/ $v\in V$ the neighborhood\/
$N(v)$ is a path or a cycle;  and in the latter case we have
 $G\cong W_4$.

\msk

\prf First, assume for a contradiction that $v$ is a vertex such
that $N(v)$ is not connected. Let $N_1$ be a component of $N(v)$,
and set $N_2=N(v)-N_1$. Let $G'$ be the graph obtained from $G$ by
deleting $v$ and introducing two new vertices $v_1$ and $v_2$
adjacent to the vertices in $N_1$ and $N_2$, respectively. In $G'$
the distance between $v_1$ and $v_2$ is at least 4. Therefore, the
inverse operation of vertex splitting can be applied to $G$, which
is a contradiction. Consequently, each set $N(v)$ induces a
connected graph.

If a vertex $u \in N(v)$ had more than two neighbors in $N(v)$, then
the edge $vu$ would belong  to more than two triangles, which
contradicts our present assumption. Hence, for every $v \in V$ and
every $u \in N(v)$, the degree  $\deg_{N(v)}(u)$ is at most two.
Since $N(v)$ is connected, this implies that $N(v)$ must be a path
or a cycle. If $N(v)\cong C_n$, then $N[v]\cong W_n$ and
$\mathcal{T}(W_n) \cong C_n$. Therefore,  $G \cong W_n$, and since
it is assumed that $G$ cannot be reduced by inverse edge splitting,
$G\cong W_4$ must be valid. \dia

\bsk

Therefore, consider only the case where $G\ncong W_4$ and $N(v)$ is
a path  for every $v\in V$. Partition the edge set of $G$ as $E = F
\cup H$ where
  $F$ is the set of doubly covered edges and
   $H$ is the set of private edges of triangles.

If $t_1,\dots,t_n$ are the triangles in the assumed cyclic order,
then we use the notation $F=\{f_1,\dots,f_n\}$ where $f_i = E(t_i)
\cap E(t_{i+1})$, and denote by $h_i$ the private edge of $t_i$.
Hence $E(t_i)=\{f_{i-1},f_i,h_i\}$ (subscript addition is
considered modulo $n$). The graphs $G_F=(V,F)$ and $G_H=(V,H)$
contain no triangles.
   It is worth noting that the original graph $G$
  can be obtained  from  $G_F$ if,
  for every
    $1\leq i\leq n$, the non-common ends of $f_i$
and $f_{i+1}$ are connected by an edge  (which is just the private
edge $h_{i+1}$ in $G$).

\bsk

 \nin
{\it Claim 2. } The graph\/ $G_F$ of doubly covered edges is a
\emph{hairy cycle} (a cycle with   any number of pendant vertices
attached to its vertices).

\msk

\prf  Let every vertex $v\in V$ be associated with the set $I(v)$ of
indices of doubly covered edges incident to $v$:
$$ i\in I(v) \qquad \Longleftrightarrow  \qquad f_i \mbox{  is incident
to }  v.$$

By definition, every index $1\leq i \leq n$ is contained in exactly
two sets $I(v)$. Since $F$ contains exactly two edges from each
triangle, every vertex of $G$ is incident with at least one edge of
$F$ and hence, no $I(v)$ is empty. The fact that  $f_i$ and
$f_{i+1}$ share a vertex  for every $1\leq i \leq n$, implies that
for every $i$ there exists a vertex $v$ such that
$\{i,i+1\}\subseteq I(v)$. The connectivity of $G_F$ also follows.

Now, consider a vertex $v\in V$, say of degree $d$. By Claim 1, its
neighborhood $N(v)$ induces a path $P_d=u_1 u_2 \dots u_d$ in $G$.
Any two consecutive vertices of $P_d$ together with $v$ form a
triangle. Hence, the doubly covered edges incident to $v$ are
exactly
 $vu_2,vu_3,\dots, vu_{d-1}$, and the $d-1$ triangles just mentioned
 correspond to consecutive vertices in the cycle $\T(G)\cong C_n$.
It also follows that the index set $I(v)$ of $v$
 contains $d-2$ consecutive integers (viewing 1 as the successor of $n$).

By these facts we obtain that the set  of vertices  which have at
least two incident doubly covered edges   induce a cycle in $G_F$.
This cycle will be referred to as $C^*=v_1v_2\dots v_k$, where the
vertices are indexed   according to the cyclic order. It contains
exactly those vertices $v_i$ for which $|I(v_i)|\geq 2$ and the two
edges of $C^*$ incident to $v_i$
 are exactly the doubly covered edges with smallest and largest
 indices from $I(v_i)$ (where `smallest' and `largest' are meant
 along a  fixed cyclic order of $1,2,\dots ,n$).

If a vertex $v_i$ from $C^*$ is incident to an edge $f_\ell=v_iu$
which does not belong to $C^*$, then $u$ cannot be incident to any
other edge of $G_F$, since both $f_{\ell-1}$ and $f_{\ell+1}$ are
also incident to $v_i$, and vertex $u$, too, must be incident to
edges with consecutive indices without a gap. Thus, all vertices and
edges not contained in the cycle $C^*$ are pendant vertices and
edges in $G_F$. Consequently, $G_F$ is a hairy cycle. \dia

\bsk

From now on,  when the inverse operation of edge splitting is
applied in $G$ to a vertex $w$ which is pendant in $G_F$, we say
that $w$ is \emph{eliminated}.
 Let us emphasize that we excluded this situation by assumption
at the very beginning; hence, several proofs below will apply
 the fact that it is impossible to identify a vertex
 which can be eliminated.

\bsk

 \nin
{\it Claim 3. } In\/ $G_F$  each vertex is
 incident with at most one pendant edge.

 \msk

  \prf Assume for a contradiction that a vertex $v_i$ has at least two
   pendant neighbors.
  Then, there exist two pendant neighbors $u$ and $w$ such that
  $f_j=v_iu$ and $f_{j+1}=v_iw$ for some $j$. Also, $f_{j-1}=v_iv'$ and
  $f_{j+2}=v_iv''$ have to be incident to $v_i$ (otherwise $u$ or $w$ would
  be a vertex from the cycle). Under these assumptions vertex $w$
  could be eliminated, since $v''w$ and $uw$ are private edges,
  moreover the only possible common neighbor, $x\neq v_i$, of $v''$ and $u$ in $G$
  might be $v'$, but since $G\not\cong W_4$, it cannot be.
    This contradiction proves the claim. \dia

\bsk

 \nin
{\it Claim 4. } If there exists a vertex\/ $v_i$ incident to a
pendant edge\/ $f_j= v_iu$ in\/ $G_F$ then the length $k$ of cycle\/
$C^*$ of\/ $G_F$ is at most 6.

\msk

  \prf
The graph $G_F$ contains no triangle, hence $v_{i-1}v_{i+1}\not\in
F$.  Moreover, since
  $f_{j-1}=v_{i-1}v_i$ and $f_{j+1}=v_iv_{i+1}$ are not consecutive
   doubly covered edges,  $v_{i-1}v_{i+1}$ cannot be a private edge of $G$ and so,
   $v_{i-1}v_{i+1}\not\in E$.
   If the only common neighbor of
 $v_{i-1}$ and $v_{i+1}$ were $v_i$, then $u$ could be eliminated and
  replaced by the edge $v_{i-1}v_{i+1}$.
 As $u$ cannot be eliminated, vertices $v_{i-1}$ and $v_{i+1}$
 have some common neighbor $x\neq v_i$ in $G$. Now, assume that
 the cycle of doubly covered edges is of
 length $k >6$. Then, every vertex of $C^*$
  is at distance greater
 than two apart from at least one of $v_{i-1}$ and $v_{i+1}$ in
 $G_F$, and the
   same is true for the possible  pendant vertices of $G_F$.
  Thus, no doubly
 covered edges and no private edges   form a second triangle with
 $v_{i-1}v_{i+1}$ in $G$, and  $u$ can be eliminated if $k>6$, which
 contradicts the present conditions. Thus, $k \leq 6$ follows. \dia

 \bsk

 If there exist no pendant edges in $G_F$ and $k\geq 7$, then $G$
 belongs to the elementary type $(b)$. Assume that this is not the case
 and $k \leq 6$.
Since $G_F$ contains no triangle, the length $k$ of $C^*$ is either
4 or 5 or 6.

\bsk

 \nin
{\it Claim 5. } If\/ $k= 4$ then\/ $G \cong S_A$.

\msk

\prf  No chord of the four-cycle belongs to the graph $G$, since it
would be a doubly covered edge and $G_F$ cannot contain triangles.
Avoiding such a case, $G_F$ has to contain precisely one pendant
edge on each vertex of $C_4$. Supplementing $G_F$ with the private
edges between the non-common ends of consecutive double edges, the
graph $G\cong S_A$  is obtained. \dia

\bsk

 \nin
{\it Claim 6. } If\/ $k= 5$ then\/ $G \cong S_B$.

\msk

\prf Let  the five-cycle be $v_1v_2v_3v_4v_5$. If one of these
vertices, say $v_i$ has no pendant edge in $G_F$, then
$v_{i-1}v_{i+1} \in H$. This edge cannot belong to any other
triangle in $G$; hence, neither $v_{i-2}v_{i+1}$ nor
$v_{i-1}v_{i+2}$ can be an edge in $G$. This means that both
$v_{i+2}$ and $v_{i-2}$ must have pendant edges in $G_F$. Similarly,
if there is no pendant edge on $v_{i+1}$, then there is one on
vertex $v_{i-1}$. This proves that if $G$ is complying with our
assumption, there exist three consecutive vertices on the five-cycle
of $G_F$ such that each of them has a pendant edge.

On the other hand, if four consecutive vertices, say $v_{i}$,
$v_{i+1}$, $v_{i+2}$ and $v_{i+3}$ have pendant edges, the pendant
vertex adjacent to $v_{i}$ can be eliminated. Consequently, if $G_F$
has a cycle of length 5, then exactly three of its consecutive
vertices have pendant edges, and $G\cong S_B$ holds.  \dia

\bsk

 \nin
{\it Claim 7. } If\/ $k= 6$ then\/ $G \cong S_C$ or\/ $G
\cong S_D$.

\msk

\prf  Consider a six-cycle   $v_1v_2v_3v_4v_5v_6$. If none of the
vertices $v_1$, $v_3$ and $v_5$ (or $v_2$, $v_4$ and $v_6$) has a
pendant edge, then a forbidden triangle $v_2v_4v_6$ ($v_1v_3v_5$)
would arise in $H$. Hence, two consecutive or two opposite vertices
of the six-cycle surely have pendant edges. On the other hand, if
two of the vertices $v_1$, $v_3$ and $v_5$ (or two of $v_2$, $v_4$
and $v_6$) have pendant edges, then one of the corresponding pendant
vertices can be eliminated.  Hence, under the given conditions,
either two consecutive or two opposite vertices have pendant edges
and $G\cong S_C$ or $G\cong S_D$ is obtained, respectively. \dia

\bsk

The above cases cover all possibilities, therefore the theorem
follows. \qed

\bsk

Note that $\T(G)$ is an odd hole if and only if
 \tmz \item $G$ can be
obtained from $W_4$, or from one of the four supplementary types, or
from $C_n^2$ where $n\geq 8$ is even, by a sequence of an odd number
of edge splittings and an arbitrarily number of vertex stickings; or
\item $G$ can be
 obtained   from $C_n^2$ where $n\geq 7$ is odd, by a sequence of an even
number of edge splittings and an arbitrarily number of vertex
stickings. \etmz

As a vertex sticking clearly creates a graph which is not
path-neighborhood, and then applying any further vertex stickings or
edge splittings this property does not change, we infer the
following characterization of
path-neighborhood graphs whose triangle graph is a cycle. This
gives a solution for Problem $3(a)$   posed by Laskar, Mulder and
Novick in \cite{Laskar1}.

\bc   \label{png-char}
 Let\/ $G$ be a graph  and assume that for every vertex\/ $v \in
V(G)$ the open neighborhood\/ $N(v)$ induces a path on at least two
vertices.
 Then,\/ $\mathcal{T}(G) \cong C_n$ for some\/
$n \geq 3$ if and only if
  \tmz
  \item[$(i)$] $G \cong C_k^2$ for a\/ $k \ge 7$; or
  \item[$(ii)$] $G \cong S_A$ or\/ $S_B$ or\/ $S_C$ or\/ $S_D$; or
  \item[$(iii)$] $G$ can be obtained  by a sequence of edge splittings
   from a graph which satisfies\/ $(i)$ or\/ $(ii)$.
     \etmz \ec

 \brm It is worth investigating the status of graphs\/ $C_n^2$ for\/
  $n\leq 6$.

\msk

\nin
 $(i)$ \
  The square\/ $C_4^2$ of the four-cycle is\/ $K_4$,
  and also its triangle graph is\/
  $K_4$. But if we double the two diagonals
  added to\/ $C_4$ for\/ $C_4^2$, and apply four edge
  splittings on these two pairs of parallel edges, we obtain\/
  $S_A$ from the supplementary  type\/ $(A)$. In fact, this does not
  correspond precisely to the definition of edge splitting, but we can
   view the case as if one of the edges\/ $v_1v_3$  belonged to triangle\/
   $v_1v_2v_3$, and its `twin' edge to\/ $v_1v_4v_3$.
   The two edges\/ $v_2v_4$ can be treated similarly.

\msk

\nin
 $(ii)$ \
   For the five-cycle\/ $C_5=v_1\dots v_5$, the square graph\/
   $C_5^2$ is a complete\/ $K_5$, and its triangle graph is the
   complement of the Petersen graph, hence not at all a cycle.
   But if we consider only the five triangles of the form\/
   $v_iv_{i+1}v_{i+2}$ (where\/ $1\leq i\leq 5$), and
   the remaining triangles of the form\/ $v_iv_{i+1}v_{i+3}$ are
   ``damaged" by properly chosen
   edge splittings, we obtain a graph\/ $G$ whose
   triangle graph is a cycle. If such a\/ $G$ cannot be reduced by the
   inverse of  edge splitting, then it is isomorphic to\/ $S_B$.

\msk

\nin
 $(iii)$ \
   A similar observation can be made for\/
   $n=6$, where at least two edge splittings have to be taken in\/
    $C_6^2$ to achieve that the triangle graph becomes
   a cycle. In the minimal (non-reducible) configurations
   exactly two edge splittings are needed,
   either neighboring or opposite. This
   yields a graph isomorphic either  to\/ $S_C$ or to\/ $S_D$.
    \erm

   \bsk

   %%%%%%%%%%%%%%%%%%%%%%%%%%%%%%

\section{Forbidden subgraph characterizations}

In Subsection~\ref{subsec:3-1} we prove the main result of this
section, which is a forbidden subgraph characterization of graphs
whose triangle graph does not contain an induced cycle of a given
length $n\ge 3$. This problem is in close connection with the
problem we solved in Theorem \ref{cycle}. Clearly, if the triangle
graph $\cT(G)$ is
 $C_n$-free, then $G$   contains no subgraph $H$ with
$\cT(H)\cong C_n$. On the other hand, avoiding these subgraphs is
not sufficient for   $\cT(G)$ to be $C_n$-free. For instance, none
of the subgraphs of $C_6^2$ has a triangle graph isomorphic to the
6-cycle, but $\cT(C_6^2)$ contains an induced $C_6$;   moreover, no
subgraph of $K_4$ has a triangle graph isomorphic to a cycle, but
$\cT(K_4)\cong K_4$ contains 3-cycles.

In Subsection~\ref{subsec:3-2} we establish some  immediate
consequences of the main result. This contains a forbidden subgraph
characterization of graphs whose triangle graph is a tree, a
chordal graph,  or a perfect graph. Further, in
Subsection~\ref{subsec:3-3} we determine a graph class on which
Conjecture~\ref{T-conj} holds.

\msk

Let us   recall that
  $C_n$-free triangle graph means that  no \emph{induced}
subgraph of $\cT(G)$ is isomorphic to $C_n$; whilst
 forbidden subgraphs given for $G$ are meant that
they must not occur as \emph{non-induced} subgraphs either.

\subsection{Graphs with $C_n$-free triangle graphs}
 \label{subsec:3-1}

In this subsection \emph{``minimal forbidden graph for $C_n$''} is
meant as a graph $G$ whose triangle graph contains an induced cycle
$C_n$ but this does not  hold for any  proper subgraph of $G$. If
$G$ is a minimal forbidden graph for $C_n$, the triangles
$T_1,\dots, T_n$ in $G$, which correspond to the vertices $t_1,
\dots , t_n$ inducing a   specified cycle in $\cT(G)$, are called
\emph{cycle-triangles}, while the other triangles of $G$  are called
\emph{additional triangles}. By definition, in a minimal forbidden
graph, every edge is contained in at least one cycle-triangle. The
edge $e \in E(G)$ will be called \emph{private edge}  if it  belongs
to exactly one cycle-triangle, and those contained in exactly two
cycle-triangles will be called \emph{doubly covered}. That is, when
an edge is declared to be
private or doubly covered, additional triangles are disregarded.

 \bp \label{n-n}
 Let\/ $n\ge 3$ be an integer.
  Assume that\/ $G$ is a minimal forbidden graph for\/ $C_n$.
   Then, either\/
 $G$ is isomorphic to\/ $K_5-K_3$, or else
   no edge of\/ $G$ is covered by more than two cycle-triangles and
 furthermore\/ $G$ contains exactly\/ $2n$ edges, from which\/ $n$ edges are private
 and\/ $n$ edges are doubly covered.
 \ep
 \pf Three cycle-triangles sharing a common edge correspond to a
 3-cycle in the triangle graph. Hence, if an edge of $G$ is covered
 by three cycle-triangles, then $n=3$. Moreover,
  by the minimality   assumption, $G\cong K_5-K_3$ must hold
   and $G$ has no additional triangles.

 Now, assume that $G \not\cong K_5-K_3$ and $G$ is minimal for
 $C_n$. Fixing an induced $n$-cycle in $\cT(G)$, each
 cycle-triangle in $G$ has exactly two neighboring cycle-triangles,
 and hence exactly two doubly covered edges. By minimality, each edge of $G$ belongs
 to at least one  cycle-triangle. Therefore, $E(G)$ consists of exactly $n$ private
 and $n$ doubly covered edges. \qed

 \bc \label{cor-n-n}
 Let $G$ be a graph which is  not isomorphic to $K_5-K_3$. Then, $G$
   is a minimal forbidden graph for\/ $C_n$ (for a specified\/ $n \ge 3$)
     if and only if\/ $\cT(G)$
 contains an induced\/ $n$-cycle and\/ $G$ has exactly\/ $2n$ edges.
  \ec
  \pf If $\cT(G)$ contains an induced $n$-cycle, then either $G$ or some
  proper subgraph of it must be minimal forbidden for $C_n$. Since
  each proper subgraph has   fewer than $2n$ edges, Proposition
  \ref{n-n} implies that   $G$ itself is a minimal forbidden graph.
  The other direction follows immediately from Proposition
  \ref{n-n}. \qed

 \bc There exists no graph which is minimal forbidden for both\/ $C_n$ and\/
 $C_m$ if\/ $n\neq m$.
 \ec

 The operations edge splitting and vertex sticking will be meant
   in the same way as
introduced in the previous section, but the conditions of their
applicability are relaxed % modified
  --- and indicated with the adjective `weak' --- as described next.
 Recall that throughout this section the position of an $n$-cycle
in $\cT(G)$ is assumed to be specified, in order to
 distinguish between cycle-triangles and additional triangles.
  We also emphasize that these operations cannot be applied for
 a graph where the fixed $n$-cycle
in $\cT(G)$ corresponds to three triangles incident to a common
edge. Particularly, we assume $G \ncong K_5-K_3$.

  \tmz
\item \emph{Weak edge splitting} can be applied for any private edge $e$.
 If this edge $e=uv$ belongs to the cycle-triangle $uvx$, we
 introduce a new vertex $w$ and change the edge set from $E$ to
 $E\setminus \{uv\} \cup \{uw, wv, wx\}$.
   This transforms each additional triangle (if exits) incident
with $e$ to a cycle of length $4$. The new vertex $w$ is of degree
$3$, and the two edges $uw$ and $wv$ originated from $e$
 are two incident private edges in the graph obtained.
 Particularly, if no additional triangles are incident with $e$, a weak edge
 splitting applied to $e$ is also called \emph{strong edge
 splitting}.

\item
 \emph{Weak vertex sticking} can be applied for any two vertices
at distance at least 3 apart. If this distance is at least 4, the
triangle graph remains unchanged and the operation is also a
\emph{strong vertex sticking}, and corresponds to  `vertex sticking'
 introduced in Section \ref{defterm}.
%% the previous section.
  A weak vertex sticking, when applied for
 vertices at distance 3, creates some
 new additional triangle(s), but a strong vertex sticking cannot cause change
in the triangle graph. \etmz

 For instance,  strong edge splitting cannot be applied for
$K_4$, but a weak edge splitting can be applied for any of its edges
and results in a wheel $W_4$. In $C_{10}^2$, no two vertices are at
distance $4$ or more, so strong vertex  sticking cannot be applied,
but  two opposite vertices of the cycle can be stuck in the weak
sense. In this case,  six additional triangles arise.

 These operations have their inverses in a natural way.
Before investigating the conditions of their applicability, let us describe
 their effect on minimal forbidden graphs.

\bp \label{weak-op}
 \tmz
 \item[$(i)$] If\/ $G'$ is obtained from\/ $G$ by a weak edge splitting
 (or, equivalently, if\/ $G$ is obtained from\/ $G'$ by an inverse weak
 edge splitting),
 then\/ $G'$ is   a minimal forbidden graph for\/ $C_{n+1}$ if and only if\/
 $G$ is minimal forbidden for\/ $C_n$.
 \item[$(ii)$] If\/ $G''$ is obtained from\/ $G$ by a weak vertex
 sticking
 (or, equivalently, if\/ $G$ is obtained from\/ $G''$ by an inverse weak
 vertex sticking),
 then\/ $G''$ is  a minimal forbidden graph for\/ $C_{n}$ if and only if\/
 $G$ is minimal forbidden for\/ $C_n$.
 \etmz
\ep
  \pf $(i)$ \, It is clear from the definition of weak edge splitting that
  the fixed $n$-cycle of $\cT(G)$ is transformed into an induced $(n+1)$-cycle
  of $\cT(G')$ and vice versa. Moreover, $|E(G)|=2n$ if and only if
  $|E(G')|=2n+2$. Hence, the statement
%% $(i)$
   follows by Corollary
  \ref{cor-n-n}.

  \msk

  \nin
  $(ii)$ \, If $G''$ is obtained from $G$ by a weak vertex
  sticking, then $\cT(G)$ contains an induced $n$-cycle
   if and only if $\cT(G'')$ contains an
  induced $n$-cycle. Additionally, $|E(G)|=2n$ holds if and only if
  $|E(G'')|=2n$. Similarly to the previous case, Corollary
  \ref{cor-n-n} implies the statement. \qed

\bsk

  Next we prove necessary and sufficient conditions  under which the
  inverse operations can be applied. Let us introduce the following
  notion. For a graph $G$ and for a fixed cycle in the triangle
  graph $\cT(G)$, the \emph{cycle-triangle neighborhood} $N^*(v)$ of a
  vertex $v\in V(G)$ is obtained by taking the vertices and edges of
  the cycle-triangles incident to $v$ and then removing vertex $v$
  and the incident edges.

  \bp \label{prop:inverse}
   Given a graph\/ $G$, with a fixed\/ $n$-cycle in its triangle
  graph such that each edge of\/ $G$ is contained in at least one
  cycle-triangle, the following statements hold:
  \tmz
   \item[$(i)$]
   An inverse weak edge splitting which  eliminates vertex\/ $w$
   exists if and only if\/ $w$ has degree\/ $3$ and  there are two
   neighbors\/
   $u$ and\/ $v$ of\/ $w$ such that\/ $uv\notin E(G)$.
        \item[$(ii)$]   An inverse strong edge splitting which  eliminates vertex\/ $w$
   exists if and only if\/ $w$ has degree\/ $3$, moreover  for the three neighbors\/
   $u$, $v$, and\/ $x$ of\/ $w$,   we have\/
 $uv\notin E(G)$,  and\/   $w$ is the only common neighbor of\/ $u$ and\/ $v$ besides\/ $x$.
    \item[$(iii)$] An inverse weak vertex sticking can be applied for a vertex\/ $v$ if and only if its
cycle-triangle neighborhood\/ $N^*(v)$ is disconnected.
  \item[$(iv)$] An inverse strong vertex sticking can be applied for a vertex\/ $v$ if and only if its
  neighborhood\/ $N(v)$ is disconnected.
  \etmz
  \ep
  \pf By definition, an edge splitting always creates a vertex of
  degree $3$. If this operation was applied for the  private edge $e=uv$ of the cycle triangle $xuv$, then
  in the obtained graph $G$, $uv$ is not an edge.
  Moreover, if   a strong edge splitting was applied, $uv$ is not contained
  in any triangles different from $xuv$. This proves that the conditions
  given in $(i)$ and $(ii)$ are necessary for the applicability of
  inverse weak and strong edge splittings.

  To prove sufficiency, first observe that a vertex $w$ which satisfies the
  conditions in $(i)$ and $(ii)$ does not belong to a $K_4$. It is also assumed that each edge is involved in a cycle-triangle.
  Then, since $w$ has three neighbors, it is incident to exactly two
  cycle-triangles, which  share an edge. This doubly covered edge must be $xw$, while the remaining two
  edges $uw$ and $vw$ are private edges in $xuw$ and $xvw$,
  respectively.  Then, by removing $w$
  and inserting the edge $uv$, the triangles $xuw$ and $xvw$ are
  replaced with $xuv$. Since the new edge $uv$ is the private edge of
  $xuv$, a weak edge splitting can be applied to it and $G$ is reconstructed.
  Under the conditions of $(ii)$ no additional
  triangle is incident to $uv$ and so $G$ can be reconstructed by a
  strong edge splitting. These prove that the inverse
  operations
  can be applied to $G$ under the conditions of $(i)$ and $(ii)$.

  To prove necessity in $(iii)$ and $(iv)$, consider two vertices
   $v_1$ and $v_2$ to which a weak vertex
  sticking is applied. By definition,
    $v_1$ and $v_2$  have distance at least $3$. Hence, there are no
  common vertices in  $N^*(v_1)$ and $N^*(v_2)$. Recall that this transformation does not create new cycle-triangles.
  Therefore,  sticking $v_1$ and $v_2$, the new vertex $v$ will have
  a   disconnected cycle-triangle neighborhood. If it is a strong vertex
  sticking, the distance of $v_1$ and $v_2$ is at least $4$ and we
  cannot have edges between the vertices of $N(v_1)$ and $N(v_2)$.
 This implies that $N(v)$ will be disconnected.

 For the other direction,   assume that the condition given in $(iii)$ holds. Let $v$ be deleted, and let the vertices from one component
 of $N^*(v)$ be joined to a new vertex $v_1$ while the further
 vertices from $N^*(v)$ be joined to another new vertex $v_2$. Then,
 the cycle triangles do not change (apart from the fact that $v$ is replaced with $v_1$ or
 $v_2$). We observe that $v_1$ and $v_2$ do not have a common
 neighbor, hence their distance is at least $3$. This proves that the original graph $G$ can be
 reconstructed
 by applying a weak vertex sticking to $v_1$ and $v_2$.
 If the stronger condition from $(iv)$ also holds for $G$, we will not have any
 edges between $N(v_1)$ and $N(v_2)$. Hence, the distance of $v_1$
 and $v_2$ is at least 4, and $G$ can be reconstructed by a strong
 vertex sticking. \qed
 \bsk

 Remark that an inverse weak edge splitting may create new additional
  triangles, an inverse weak vertex sticking may damage some of the
  additional triangles, while the  inverse
  strong vertex sticking keeps the triangle graph the same.

 Concerning the order in which   these transformations can be
applied, we prove the following property.
 (Although the third part
could also be made more detailed,
 by performing strong vertex stickings before non-strong ones,
 we do not need this fact in the current context.)

 \bp \label{order}
 Assume that graph\/ $F$ can be obtained from\/ $G$ by a sequence of
 weak edge splittings and weak vertex stickings. Then,\/ $F$   can also
 be obtained from\/ $G$ by performing the operations in the following
 order:

 \nmr
 \item some (maybe zero) weak but not strong edge splittings,
 \item some (maybe zero)   strong edge splittings,
 \item some (maybe zero) weak vertex stickings.
 \enmr
 Consequently,\/ $G$ can be obtained from\/ $F$ in the reverse order of
 the corresponding inverse operations.
 \ep
 \pf First, observe that if   the $i$th transformation $O_i$ is a weak sticking of
 $x$ and $y$, moreover the $(i+1)$st transformation $O_{i+1}$ is the weak splitting of
 the edge $e$, then they   can also be applied in the order $O_{i+1},  O_i$.
  Indeed, a vertex sticking does not create a new private edge,
 hence $O_{i+1}$ can be performed before $O_i$. On the other hand,
 an edge splitting cannot decrease the distance of $x$ and $y$ and
 cannot create a new edge between two vertices which were present
 previously. Hence, the order $O_{i+1},  O_i$ is feasible and
 gives the same result as $O_{i},  O_{i+1}$.
  Therefore, we can   re-order the transformations in such a way that
 all edge splittings precede all vertex stickings.

 A sequence of edge splittings  can be
  unambiguously
  described by assigning  a nonnegative
 integer $s(e)$ to each edge $e \in E(G)$, where $s(e)$ is the
 number of edge splittings applied to $e$   and to the edges originated
 from $e$. Equivalently, this is the number of
  subdivision vertices we have on $e$ at the end.
  This also shows that edge splittings can be performed in any
  order. If we want to start with weak but not strong edge
  splittings, we just take an edge with $s(e) \ge 1$ which is incident
  with an additional triangle and apply an edge splitting as long as
  such an edge exists.  \qed

 \thm \label{main}
   Let\/ $n\ge 3$ be a given integer.
 The triangle graph\/ $\cT(G)$ of a graph\/
 $G$ does not contain an induced cycle of length\/ $n$ if and only if\/
 $G$ has no subgraph which is isomorphic to   any of the following forbidden
 ones.
  \tmz
  \item[$(a)$] If\/ $n=3$, the forbidden subgraphs are\/ $K_4$ and\/ $K_5-K_3$.
  \item[$(b)$] If\/ $n=4$, the only forbidden subgraph is\/ $W_4$.
  \item[$(c)$] If\/ $n=5$, the forbidden subgraphs are\/ $W_5$ and\/ $C_5^2\cong K_5$.
 \item[$(d)$] If\/ $n=6$, the forbidden subgraphs are\/ $W_6$,\/ $C_6^2$,\/
  $K_6-K_3$,\/ $K_6-P_4$, and the graph obtained from\/ $C_5^2\cong K_5$ by a weak
  edge splitting.
  \item[$(e)$] If\/ $n\ge 7$, the forbidden subgraphs are
    \tmz
    \item[$(i)$] $W_n$;
    \item[$(ii)$] graphs   obtained from\/ $C_m^2$ by\/ $n-m$ weak edge splittings
         (\/$5 \le m \le n$);
    \item[$(iii)$] graphs   obtained from\/ $K_6-K_3$ by\/ $n-6$ weak edge splittings;
    \item[$(iv)$] graphs   obtained from\/ $K_6-P_4$ by\/ $n-6$ weak edge splittings;
    \item[$(v)$] graphs obtained from any graphs described in\/
    $(ii)-(iv)$ by any number of weak vertex stickings.
    \etmz
  \etmz
 \ethm

    \pf First, observe that   three triangles, any two of which
    share an edge, are either three triangles having a fixed
    common edge,
   or they belong to a common
    $K_4$ subgraph. Hence, any $3$-cycle in $\cT(G)$ origins  either from a
    $K_5-K_3$ or from a $K_4$.
      Consequently, by minimality, if $n=3$ then either
    $G\cong K_5-K_3$ or $G\cong K_4$ holds.
    From now on, we
    consider a graph $G$ which is not isomorphic to $K_5-K_3$ and
  is  minimal forbidden for a specified $n\ge 3$.
   By Proposition \ref{n-n}, there are exactly $n$
    private edges and $n$ doubly covered edges in $G$.

    By Proposition \ref{weak-op}, weak edge splittings, vertex stickings and their
    inverse operations do not change the status of a graph being
    minimal forbidden for at least one cycle $C_n$. Hence,  we may assume further that   inverse weak vertex sticking and
    inverse weak edge splitting cannot be applied to $G$.

      We have the following two cases concerning the additional
     triangles of $G$.

        \paragraph{\bf{Case 1.}} Each   additional triangle
     contains at least one private edge (or there is no additional triangle).

   In this case, we apply a minimum number of weak edge splittings such that   all the additional
    triangles of $G$
    are damaged.  This yields a graph $G'$ with $\cT(G')\cong C_n$.
         We shall prove that neither inverse strong edge splitting nor inverse strong
    vertex sticking can  be applied for $G'$.
    \bsk

    By our assumption, inverse weak vertex sticking cannot be applied for $G$. Hence, by Proposition~\ref{prop:inverse}$(iii)$,
    for every $v\in V(G)$, $N^*(v)$ is connected. Assume first that an edge
    splitting is applied for an edge $e=vu$ of $G$. Then, $u$ is omitted
    from the neighborhood of $v$, but the new
    vertex $w$ appears    in   $N^*(v)$ and  has exactly the same
    neighbor there as $u$ had.
    Therefore,  $N^*(v)$ remains connected. Now, consider an edge splitting applied for  a private edge
    $xy$ from the cycle-triangle $vxy$. In   $N^*(v)$, this
    means only the subdivision of  the edge $xy$. This also keeps
    connectivity. As the third case, for any new vertex $w$  which was introduced by an edge
    splitting, $N^*(w)$ is a path of order $3$.
    Therefore, every vertex of $G'$ has a connected cycle-triangle
    neighborhood, and by Proposition~\ref{prop:inverse}$(iii)$, inverse weak (and also, strong) vertex sticking cannot be applied
    for $G'$.
    \bsk

    Concerning the other operation, we supposed that inverse weak edge splitting cannot be applied for
    $G$.
    Then we applied minimum number of weak but not strong edge
    splittings to damage all the additional triangles. The minimality condition implies that
    each new vertex belongs to at least one induced $4$-cycle originated
    from an additional triangle. Thus,  every new vertex $x$ has
     two neighbors $x_1$ and $x_2$ which share a neighbor $y$ such
     that $y$ is not adjacent to $x$. By
     Proposition~\ref{prop:inverse}$(ii)$, inverse strong edge
     splitting cannot be applied for $x$.
     The second case is when   a vertex $w$ was present already in
     $G$ and had degree $3$.
      By Proposition~\ref{prop:inverse}$(i)$, as inverse weak edge splitting cannot be applied for $G$,
      the neighbors of $w$ are pairwise adjacent. Let $u$ and $v$ be the neighbors of $w$ such that
      $uw$ and $vw$ are the private edges of triangles $xuw$ and $xvw$, respectively.
             While  minimum number of edge splittings were performed,
       no edge could be split twice. Hence, if  weak edge splitting was applied for  neither $uw$ nor
       $vw$, then $u$ and $v$ either remain adjacent or
       have a common neighbor in $G$ that is different from $x$ and $w$ (the latter case
       occurs when the edge $uv$ was split).  Then, Proposition~\ref{prop:inverse}$(ii)$ implies that vertex $w$ cannot be eliminated by an inverse
       strong  edge splitting.
   Now, assume that at least one of
         the edges $uw$ and $vw$, say $uw$ was split by inserting a new edge $xu'$.
          By our minimality condition, the
         new vertex  $u'$ is contained in an induced $4$-cycle.
         Since $u'$ has only three neighbors $u$, $x$, and $w$, furthermore $xu, xw \in E(G')$,
         the induced $4$-cycle contains $u$, $u'$ and $w$ plus one
         vertex  which is different from
         $x$. This fourth vertex must be $v$, because $w$ is also of degree $3$.
         Again, the two neighbors of $w$, namely  $u'$ and $v$, have the common neighbor which is  $u$.
           Thus, by Proposition~\ref{prop:inverse}$(ii)$, $w$ cannot be eliminated by an inverse strong edge splitting.
       Finally, we observe that the edge  splittings performed in $G$
       do not decrease the degrees of the vertices. Hence, if
        a vertex has
       degree greater than 3 in $G$, it cannot be eliminated by an inverse strong edge
       splitting in $G'$.

       Therefore, inverse strong edge splitting and inverse strong vertex sticking cannot be applied for
       $G'$. By
       Theorem~\ref{cycle}, graph $G'$ is  isomorphic
     either to $W_4$, or to $C_n^2$ with  $n\ge 7$, or to one of
     the supplementary types $S_A$, $S_B$, $S_C$, $S_D$. According to the way $G'$ is derived,
       we see that either $G=G'$ or
     $G$ can be reconstructed from $G'$ by applying some number of inverse weak edge
     splittings, to be performed as long as at least one is possible, because
      it has been assumed that inverse weak edge splitting cannot be applied to $G$.

     Checking all items from our list for $G'$, we can observe the following.
     \tmz
     \item In $W_4$, we can apply inverse weak edge splitting exactly once. This yields  $K_4 \cong W_3$.
     \item In $C_n^2$ (with $n\ge 7$) there are no vertices of degree $3$. Thus, inverse weak edge splitting cannot be
     applied.
     \item In $S_A$ we have four vertices of degree $3$, and we can choose from four possible inverse weak edge
     splittings at the first step (all the four are isomorphic). After one is performed, only
     two further (isomorphic) possibilities remain. At the end, after two inverse edge splittings we
     obtain $K_6 - P_4$.
     \item In $S_B$ three vertices have degree $3$, and  inverse weak edge splittings can be
     applied to all of them. These can be performed in any order, the result will be
      $C_5^2 \cong K_5$.
      \item For $S_C$ and $S_D$ we can apply inverse weak edge
     splitting twice. After performing them we have a $C_6^2$.
     \etmz
 %%    As we assumed, no inverse weak edge splittings and vertex
 %%    stickings can be applied for $G$. Hence,
     We conclude that in this case $G$ must
     be isomorphic either to $K_4$, or to $C_n^2$ with $n\ge 5$, or to $K_6 - P_4$.

    \paragraph{\bf{Case 2.}} There is at least one
    additional triangle $uvw$ in $G$ such that each of  the edges
    $uv$, $vw$, $uw$ is doubly covered.

    By our assumption,
    no inverse weak vertex sticking  can be applied for $G$. Hence,
    every vertex $x$ has a connected cycle-triangle neighborhood. Since every
    triangle $T_i$ has exactly one private edge and the two
     doubly covered edges correspond  to $T_{i-1} \cap T_i$ and $T_i
     \cap T_{i+1}$, the triangles incident with $x$ are consecutive
     triangles   along the triangle cycle.
%% (in cyclic order).
      Let us refer to this property as `continuity'.

     Since $uv$, $vw$ and  $uw$ are doubly covered, the incident
     triangles can be given with their vertex sets in the form  $T_1=
     \{uva_1\}$, $T_2=\{uva_2\}$, $T_i=\{vwb_1\}$,
     $T_{i+1}=\{vwb_2\}$, $T_j=\{wuc_1\}$ and $T_{j+1}=\{wuc_2\}$,
     where   these six triangles are not assumed to be consecutive, but
     they are given in a cyclic order.

     If $a_2=b_1$, $b_2=c_1$ and $c_2=a_1$, we get  the desired
     result $G\cong K_6-K_3$. Now, suppose that $a_2\neq b_1$,
     which is equivalent to $i \neq 3$. Since $T_2$ and $T_i$   are
     incident with vertex $v$ (but  $T_j$ is not),   the
     continuity of triangles at $v$ implies that $T_3$ has vertices
     $va_2x$, moreover $a_2x$ must be the private edge. Now, $ua_2$
     and $a_2x$ are two incident private edges of consecutive
     triangles. The inverse weak edge splitting could not be
     applied for them only if $ux$ is an edge in $G$.
     Then,  the triangle $T_\ell$ having this edge $ux$ belongs to
     triangles incident with $u$. Since $T_1,T_2$ are incident with
     $u $ but $T_3, T_i, T_{i+1}$ are not, $i+2 \le \ell$ follows.
     But in this case, the triangles incident with $x$ cannot
     satisfy the continuity (since  $x$ belongs to $T_3$, might belong to $T_i$, but it is
     surely not contained in $T_{i+1}$ and $T_1$).
     This contradiction proves that $a_2=b_1$ and similarly, $b_2=c_1$ and
     $c_2=a_1$ must be   valid, as well. Thus, $G\cong K_6-K_3$ holds.

\msk

     These cases together cover all possibilities, therefore the theorem is
     proved. \qed

\subsection{Trees, chordal graphs and perfect graphs}
  \label{subsec:3-2}

\bl $\T(G)$ is connected if and only if there does not exist a
partition of\/ $E(G)$ into two sets\/ $A$ and\/ $B$ such that  each
of\/ $A$ and\/ $B$ contains at least three edges which   induce a
triangle and each triangle in\/ $G$ is either in\/ $A$ or in\/ $B$.
\el

\pf   Let $A\cup B=E(G)$ be an edge partition such that
 each triangle of $G$ is contained in either $A$ or $B$.
Then in $\cT(G)$ there cannot be any edges from the vertices representing
 the triangles inside $A$ to those representing the triangles inside $B$.
Thus, if there exist two triangles $T_A\subset A$ and $T_B\subset B$ in $G$,
 then $\T(G)$ has at least two components.

Conversely, assume that $\T(G)$ is disconnected. Let $A$ be the
collection of all edges of $G$ corresponding to the triangles in
one of the components
%% $C$
 of $\T(G)$, and let $B = E(G) - A$. This $B$
also contains at least one triangle, since $\T(G)$ is
disconnected. Now, $\{A,B\}$ is a partition of $E(G)$ such that
each of $A$ and $B$ contains at least one triangle and each triangle
is either in $A$ or in $B$. \qed

\bsk

  Let us say that graph $G$ is {\it triangle-connected\/} if $\T(G)$ is connected.

\bsk

Now, the   characterization of
 graphs whose triangle graph is a tree or a chordal graph
follows immediately from Theorem \ref{main}.

\bc For a graph\/ $G$, its triangle graph\/ $\cT(G)$ is a tree if and
only if\/ $G$ is triangle-connected and does not contain a subgraph
which is isomorphic to one of the following graphs.
    \tmz
    \item[$(a)$] $W_n$, for\/ $n\ge 3$;
    \item[$(b)$] $K_5-K_3$;
    \item[$(c)$] $C_n^2$, for\/ $n\ge 5$;
    \item[$(d)$] $K_6-K_3$;
    \item[$(e)$]  $K_6-P_4$;
    \item[$(f)$] graphs obtained from any of the graphs described in\/
     $(c)-(e)$ by any   number of weak edge splittings and weak vertex
    stickings.
    \etmz
\ec

\bc For a graph\/ $G$, its triangle graph\/ $\cT(G)$ is chordal if
and only if\/ $G$ does not contain a subgraph which is isomorphic to
any of the following graphs:
    \tmz
    \item[$(a)$] $W_n$, for\/ $n\ge 4$;
    \item[$(b)$] $C_n^2$, for\/ $n\ge 5$;
    \item[$(c)$] $K_6-K_3$;
    \item[$(d)$]  $K_6-P_4$;
    \item[$(e)$] graphs obtained from any graphs described in\/ $(b)-(d)$
     by any   number of weak edge splittings and weak vertex
    stickings.
    \etmz
\ec

Imposing parity conditions, we also obtain a characterization of graphs
 whose triangle graph is perfect.

\thm \label{perfect} For a graph\/ $G$, its triangle graph\/
$\cT(G)$ is perfect if and only if\/ $G$ does not contain   any
subgraph which is isomorphic to one of the following graphs:
    \tmz
    \item[$(a)$] $W_n$, for an odd integer\/ $n\ge 5$;
    \item[$(b)$] graphs obtained from\/ $C_n^2$ by   an even number
%% (at least 0)
     of weak edge splittings for an odd\/ $n\ge 5$;
    \item[$(c)$] graphs obtained from\/ $C_n^2$ by an odd number
       of weak edge splittings for an even\/ $n\ge 6$;
    \item[$(d)$] graphs obtained from\/ $K_6-K_3$ by an odd number of weak edge splittings;
    \item[$(e)$] graphs obtained from\/ $K_6-P_4$ by an odd number of weak edge splittings;
    \item[$(f)$] graphs obtained from the graphs described in\/ $(b)-(e)$
     by any number of weak vertex
    stickings.
    \etmz
\ethm

\pf Since $\overline{K_2} \vee \overline{P_3}$
  is forbidden for triangle
graphs \cite{Eri} and $\overline{C_n}$ contains it as an induced
subgraph for all $n \ge 7$, we have that $\T(G)$ is
$\overline{C_n}$-free for $n\ge7$. Also, $\overline{C_5}=C_5$.
Therefore, by the Strong Perfect Graph Theorem \cite{Chu}, $\T(G)$
is perfect if and only if has no induced odd hole. Moreover, $\T(G)$
contains an induced odd hole if and only if $G$   has a  subgraph
from the types described in $(a)-(f)$. This completes the proof.
\qed

\subsection{Consequences for triangle packing and covering}
  \label{subsec:3-3}

 Here we consider Conjecture~\ref{T-conj} which was posed
% by Tuza
  in \cite{Tuz81}.
 To discuss it in a more detailed way,  we need some definitions.
  We say that a family $\mathcal{F}$ of triangles in $G=(V,E)$ is {\it
independent\/} if the members of $\mathcal{F}$ are pairwise
edge-disjoint. An edge set $E'\subseteq E$ is a {\it
$\cT$-transversal\/} if every triangle of $G$ contains at least one
edge from $E'$. We denote by $\nut(G)$ the maximum cardinality of an
independent family of triangles in $G$, and by $\taut(G)$ the
minimum cardinality of a $\cT$-transversal in $G$. With this notation,
Conjecture~\ref{T-conj} is equivalent to the statement
$\tau_\Delta(G) \le 2\nu_\Delta (G)$.

This inequality has been proved only for few classes of graphs; namely, for
planar graphs, some subclasses of chordal graphs,
 graphs with $n$ vertices and at least   $\frac{7}{16}n^2$ edges
\cite{Tuz2}, graphs without a subgraph homeomorphic to $K_{3,3}$
 \cite{Kri}, graphs with chromatic number three \cite{HK},
 graphs in which every subgraph has average degree smaller
 than seven \cite{Pue},
 odd-wheel-free graphs, and graphs admitting an edge 3-coloring
 in which each triangle receives three distinct colors on its
 edges~\cite{ABT}. (The latter class contains all graphs with
 chromatic number at most four, and also all graphs which have
 a homomorphism into   the third power  of an even cycle,
  $C^3_{2k}$ with $k\ge 5$.)

The case of equality $\tau_\Delta=\nu_\Delta$ has also been
 studied to some extent (\cite{T94,ABT}).
For instance, it was proved in \cite{ABT} that
 $\tau_\Delta(G)=\nu_\Delta (G)$ is valid for
 $K_4$-free graphs $G$  whose triangle graph is odd-hole free.
Now, by Theorem \ref{perfect}, this class of graphs is
determined in a more direct way by forbidden subgraph
characterization. (Some redundancies could be eliminated; e.g.,
$K_6-P_4$ and $K_6-K_3$ contain $K_4$ which is forbidden, too.
But after some appropriate edge splittings, the $K_4$ subgraphs
disappear, hence the graphs listed in parts $(d)$ and $(e)$ of Theorem
\ref{perfect} cannot be totally omitted.)

Here we prove Conjecture~\ref{T-conj} for graphs $G$ whose triangle
graph is perfect, but the $K_4$ subgraphs are not excluded from $G$.
The condition that $\cT(G)$ is perfect, can be replaced either by
assuming that $\T(G)$ is odd-hole-free
 or by the forbidden
subgraph characterization of Theorem \ref{perfect}.

\thm   \label{perf-taunu}
 If the triangle graph of a graph\/ $G$ is perfect, then\/
$\tau_\Delta(G) \le 2\nu_\Delta (G)$ holds.
 \ethm

\pf As discussed already in \cite{ABT}, the maximum number
$\nu_\Delta (G)$ of independent triangles in $G$ equals the
independence number $\alpha(\cT(G))$ of the triangle graph. Since
the triangle graph is supposed to be perfect, its complement is also
perfect and we have
$$\nu_\Delta (G)=\alpha(\cT(G))=
\omega(\overline{\cT(G)})=\chi(\overline{\cT(G)})=\theta(\cT(G)),$$
where $\theta(\cT(G))$ is the minimum number of cliques in the
triangle graph which together cover all vertices.

In a triangle graph $\cT(G)$ we may have two types of cliques: $(A)$
its vertices correspond to triangles of $G$ all of which are
incident with a fixed edge; $(B)$ its vertices correspond to  four
triangles of a $K_4$ subgraph of $G$.

Having a minimum clique cover of  $\cT(G)$ at hand we can construct
an edge
 cover for triangles of $G$. For every clique $C_A$ of type $A$, we put the
 corresponding edge of $G$ into the covering set. This edge covers all triangles
 corresponding to the vertices covered by $C_A$ in $\cT(G)$. Then, for every
 clique $C_B$
 of type $B$, we put two independent edges from the corresponding
 $K_4$ subgraph of $G$ into the covering set. These two edges together
 cover all the four triangles.
 Since every vertex $t\in V(\cT(G))$ is covered by a clique in
 $\cT(G)$, every triangle of $G$ is covered by at least one of the selected edges.
 Thus, the set of the selected edges (wchich covers all triangles of
 $G$) contains at most $2\theta(\cT(G))=2\nu_\Delta (G)$
 edges. Consequently, $\tau_\Delta(G) \le 2\nu_\Delta (G)$ is valid.
 \qed

\paragraph{Note added on 4$^{\mbox{\scriptsize\bf th}}$ November, 2014.}

After the first appearance of this manuscript,
 Gregory Puleo kindly informed us that from his results in \cite{Pue}
 the inequality $\tau_\Delta(G) \le 2\nu_\Delta (G)$ follows for a class of graphs
  which is larger than the one in our Theorem \ref{perf-taunu}.
Namely, more generally than the graphs with perfect $\cT(G)$,
 it suffices to assume that $G$ has no subgraph isomorphic to
 $W_n$ for any odd $n \ge 5$.
That is, from the forbidden subgraphs listed in
 our Theorem \ref{perfect}, already the case (a) is sufficient to derive
 $\tau_\Delta \le 2\nu_\Delta$.
As Puleo explains in his email, this follows by the properties of the so-called
 `weak K\"onig--Egerv\'ary graphs',
%% applying Proposition 4.8 and Lemma 4.4
  proved in Section 4 of \cite{Pue}.
%%, in which no assumption is put on vertex degrees.
This extends Theorem 3 of \cite{ABT} where the analogous result was proved
 for graphs without any odd wheels (i.e., excluding $W_3\cong K_4$, too).

\paragraph{Acknowledgements.}

The authors thank the referees for their comments, in particular for
 proposing to extend a less detailed argument for Case 1 in the proof of
 Theorem \ref{main}.
Research of the first author was supported in part by
 University Grants Commission, India  under the grant MRP(S)-0843/13-14/KLMG043/UGC-SWRO.
Research of the second and third authors was
 supported in part
   by the Hungarian NKFIH
   grant SNN 1160959, and
   by the Hungarian State and the European Union
   under the grant
 T\'AMOP-4.2.2.A-11/1/KONV-2012-0072.
\bsk

\end{document}